\documentclass[12pt]{amsart}
\usepackage{latexcad}
\newtheorem{thm}{Theorem}[section]
\newtheorem{cor}[thm]{Corollary}
\newtheorem{lem}[thm]{Lemma}
\newtheorem{prop}[thm]{Proposition}
\theoremstyle{definition}

\theoremstyle{remark}
\newtheorem{rem}[thm]{Remark}

\numberwithin{equation}{section}

\newcommand{\CC}{\mathbb C}
\newcommand{\LL}{\mathbb L}
\newcommand{\NN}{\mathbb N}
\newcommand{\RR}{\mathbb R}

\newcommand{\re}{\hbox{\rm Re }}

\newcommand{\cc}{{\mathbf c}}
\newcommand{\vv}{{\mathbf v}}
\newcommand{\xx}{{\mathbf x}}

\newcommand{\cR}{\mathcal{R}}
\begin{document}

\title{CR manifolds admitting a CR Contraction}%
\author[K.-T. Kim and J.-C. Yoccoz]
{Kang-Tae Kim and Jean-Christophe Yoccoz}%

\address{(Kim) Department of Mathematics, Pohang University of
Science and Technology, Pohang 790-784 The Republic of Korea}%
\email{kimkt@postech.ac.kr}%
\address{(Yoccoz) Coll\`ege de France, 3 rue d'Ulm, 75005 Paris,
France}%
\email{jean-c.yoccoz@college-de-france.fr}%

\thanks{The research of the first named author has been supported
in part by the Grant R01-2005-000-10771-0 from the Korea Science
and Engineering Foundation.}%

\begin{abstract}
We classify the germs of $\mathcal{C}^\infty$ CR manifolds that
admit a smooth CR contraction. We show that such a CR manifold is
embedded into $\CC^n$ as a real hypersurface defined by a
polynomial defining function consisting of monomials whose degrees
are completely determined by the extended resonances of the
contraction. Furthermore the contraction map extends to a
holomorphic contraction that coincides in fact with its polynomial
normal form.  Consequently, several results concerning Complex and
CR geometry are derived.
\end{abstract}
\maketitle

\section{Introduction}

\subsection{}
A {\it CR manifold} is defined to be a smooth ($\mathcal{C}^\infty$)
manifold, say $M$, equipped with a CR structure, i.e., a sub-bundle
$\LL$ of its complexified tangent bundle $T^{\CC} M$ satisfying two conditions:
\begin{itemize}
\item[($\romannumeral1$)] $\LL_p \cap \overline{\LL}_p = \{0\}$, and %
\item[($\romannumeral2$)] $[L_1, L_2]$ is a section of $\LL$ whenever
$L_1$ and $L_2$ are.
\end{itemize}
Such a CR structure is in fact called {\it integrable}. But we restrict
ourselves in this paper to $\mathcal{C}^\infty$ smooth manifolds
with an integrable CR structure only.
\medskip

The CR structure defines a distribution $(p \in M) \mapsto \zeta_p:=
\re (\LL_p \oplus \overline{\LL}_p)$ with a {\it complex
structure} that is a smooth assignment $p \in M \mapsto J_p:
\zeta_p \to GL(\zeta_p)$ satisfying $J_p \circ J_p = -id$,
where $id$ represents the identity map. In fact,
$J_p$ for each $p$ is the almost complex structure on $\zeta_p$
induced by the action of $\CC$ on $\LL$. It is known that the CR
structure $\LL$ is uniquely determined by $J$.
\medskip

Typical CR manifolds came originally from the smooth boundaries of
a domain in $\CC^n$.  In such a case the boundary has real
dimension $2n-1$ and the distribution $\zeta$ consists of $2n-2$
dimensional hyperplanes. Accordingly a CR manifold is said to be
of {\it hypersurface type} if $\dim_\RR M = 2n-1$ and $\dim_\RR
\zeta_p = 2n-2$ for every $p \in M$.  Again, we restrict ourselves
only to the CR manifolds of hypersurface type.  Thus from here on,
by a CR manifold we always mean a CR manifold of hypersurface
type.
\medskip

A {\it CR map} is a smooth mapping, say $f$, of a CR manifold
$(M,\zeta,J)$ into another $(N, \xi, \widetilde J)$ satisfying $df
\circ J = \widetilde J \circ df$.
\medskip

\subsection{}
The purpose of this article is to study the germ $(M,p)$ of a
$\mathcal{C}^\infty$ smooth CR manifold $M$ at a point, say $p$,
that admits a smooth CR contraction $f:(M,p) \to (M,p)$. Here, a
CR map $f$ is said to be a {\it contraction}, if every eigenvalue
of $df_p:T^\CC_p M \to T^\CC_p M$ has modulus strictly smaller
than 1.
\medskip

For a germ of a CR manifold, it is quite a strong assumption to
admit a CR contraction, and yet there is a reasonably broad
collection of examples. Hence it is natural to attempt to classify
them; that is exactly the purpose of this article.
\medskip

The contractions admit a polynomial normal form (up to a
conjugation by a smooth diffeomorphism), as is well-known,
depending upon the resonances between the eigenvalues. The main
results of this paper (cf.\ Theorem \ref{theorem-main-2} for
precise presentation) are as follows:
\begin{itemize}
\item[($\romannumeral1$)] The CR manifold germ admitting a CR
contraction is CR equivalent to a CR hypersurface germ in $\CC^n$
(consequently, it is embedded) and the CR contraction extends to a
holomorphic contraction.
\medskip

\item[($\romannumeral2$)] The CR equivalence, mentioned in
($\romannumeral1$), conjugates the CR contraction to a polynomial
mapping that is in fact its normal form.
\medskip

\item[($\romannumeral3$)] The defining function of the CR
hypersurface germ (the image of the embedding mentioned in
($\romannumeral1$)) is a weighted homogeneous real-valued
polynomial each of whose monomial terms has degrees completely
determined by the {\it extended resonances} (see Section 2) of the
contraction.
\end{itemize}

\subsection{Organization of this paper}
We present the contents in the following order.
\medskip

First we recall several known facts concerning the contractions.
Some statements have been modified so that it can fit to the
purpose of this paper.  However, we point out that there is no
novelty by us in this part; we collect them and compose them as an
appendix (See Section 5) to this paper for the conveniences of the
readers.
\medskip

Second, we analyze (in Section 2) the real hypersurface germ in
$\CC^n$ at the origin that is invariant under a holomorphic
contraction.
\medskip

Third, we show (in Section 3) that any abstract smooth CR manifold
germ admitting a CR contraction is CR embedded as a hypersurface
in $\CC^n$. In fact we establish the aforementioned main results
(i), (ii) and (iii).
\medskip

Section 4 then lists some of the applications. It may be worth
noting that in contrast to the preceding related achievements,
this paper takes the existence of a CR contraction as the
governing environment. Then we derive the complex and CR geometric
consequences. It appeals to us that this change of viewpoint is
worthwhile; in particular, our classification of CR manifolds
includes even the ones that contain complex lines.

\section{Real hypersurface germ in $\CC^n$ invariant under a
holomorphic contraction}

\subsection{A holomorphic contraction}
Let $f: (\CC^n, 0) \to (\CC^n, 0)$ be a local germ of a
holomorphic contraction with a fixed point at the origin $0$.
Denote by $\lambda_1, \ldots, \lambda_m$ the {\sl distinct}
eigenvalues of the derivative $Df_0 $ of $f$ at $0$ indexed in
such a way that
$$
|\lambda_m| \le |\lambda_{m-1}| \le \cdots \le |\lambda_1|<1.
$$
Let $\nu \in \{1, \ldots, m\}$, and let $E_\nu$ represent the
characteristic subspace corresponding to $\lambda_\nu$. Then it
holds obviously that $\CC^n = E_1 \oplus \cdots \oplus E_m$. We
shall now use a vector coordinate system $\vv_1, \ldots, \vv_m$
for $\CC^n$ so that each $\vv_\nu$ gives a vector coordinate
system for $E_\nu$ for every $\nu$.

\subsection{Resonances and Extended Resonances}
Let $\NN$ denote the set of positive integers.
For each $\nu \in \{1, \ldots, m\}$ consider {\it the set of
resonances}
$$
\cR_\nu := \{I = (I_1, \ldots, I_m) \in \NN^m \mid \sum I_\ell \ge
2, \lambda_\nu = \lambda^{I_k}, \forall k=1,\ldots,m \}.
$$
\medskip

For the purpose of this article, one needs another type of
resonances.  Let
\begin{multline*}
\cR'_\nu = \{ (I,I') = (I_1, \ldots, I_m; I'_1, \ldots, I'_m)
\in \NN^m \times \NN^m \mid \\
\sum I_\ell + \sum I'_\ell \ge 2, \lambda_\nu = \lambda^I
\overline{\lambda^{I'}} \}
\end{multline*}
and call it {\it the set of extended resonances}.
\medskip

\begin{rem} \rm
(1) The resonance set $\cR_j$ and the generalized resonance set
$\cR_j'$ are finite. \medskip

(2) If $J \in \cR_j$, then $J_k = 0$ for $k \ge j$.  Similarly if
$(J,J')$ belongs to $\cR_j'$, then $J_k = 0 = J_k'$ for $k \ge j$.
\end{rem}
\medskip

\subsection{Normalization of the Contraction $f$}
By the basic theorem on normal forms for contractions (see
Proposition \ref{contraction} in \S\ref{contraction-section}),
after a holomorphic change of variables, one can write:
$$
f(\vv_1, \ldots, \vv_m) = (\widehat \vv_1, \ldots, \widehat \vv_m)
$$
with
$$
\widehat \vv_\nu = \Lambda_\nu \vv_\nu + \sum_{J \in R_\nu} \cc_J
(\vv),
$$
where:
\begin{itemize}
\item $\Lambda_\nu :E_\nu \to E_\nu$ is linear with a single
eigenvalue $\lambda_\nu$
\item $\cc_J:\CC^n \to E_\nu$ is a polynomial with values in $E_\nu$,
homogeneous of degree $J_k$ in the variable $\vv_k$, for each $k
\in \{1, \ldots, m\}$.
\end{itemize}
\bigskip

\subsection{$Df_0$-invariant hyperplanes}
Let $(M,0)$ represent the germ of smooth real hypersurface in $\CC^n$
and $f$ the holomorphic contraction above.  We work now with the
assumption that $f((M,0)) \subset (M,0)$ throughout this section.
\medskip

The following lemma is immediate:
\bigskip

\begin{lem} \label{tangent-lemma} \sl
Consider the action of the transposed operator $(Df_0)^*$ on the
real-dual $(\RR^{2n})^* = E_1^* \oplus \cdots \oplus E_m^*$. If a
non-zero dual vector $\varphi \in (\RR^{2n})^*$ satisfies that
$(Df_0)^* (\varphi) = \lambda \varphi$ for some real number
$\lambda$, then there exists an index $i \in \{1, \ldots, m\}$
such that $\lambda = \lambda_i$ and $\varphi \in E_i^*$.
\end{lem}

Since $(M,0)$ is the germ of a smooth real hypersurface through
the origin that is invariant under $f$, the (extrinsic) tangent
space $T_0 M$ is invariant under the action of $Df_0$. Let
$\varphi$ be a non-vanishing real linear-form whose kernel is $T_0
M$. Then we must have $(Df_0)^* (\varphi) = \lambda \varphi$ with
the constant $\lambda$ real.

By Lemma~\ref{tangent-lemma}, there exists $i \in \{1, \ldots,
m\}$ with $\lambda = \lambda_i$ and a decomposition $E_i = \CC e_i
\oplus E_i'$ (with coordinates $\vv_i = (z_i, {\vv'}_i)$) of $E_i$
such that
\begin{itemize}
\item $T_0 M$ is the real hyperplane defined by $\text{Re } z_i =
0$.
\item $\Lambda_i \vv_i = (\lambda_i z_i, {\Lambda'}_i {\vv'}_i +
z_i \cc')$ for some $\cc' \in E'_i$.
\end{itemize}
\bigskip

Here, we write as: $z_i = x_i + \sqrt{-1}y_i$ with $x_i, y_i \in
\RR$, and $\cc_I = (c_I^0, {\cc'}_I)$ with $I \in R_i$.  \sf This
special index $i$ will be kept separately from the other indices
here as well as in what follows. \rm
\bigskip

Now we are ready to present the main result of this section.

\begin{thm} \label{theorem-main}
The germ $(M,0)$ of a smooth real hypersurface at $0$ invariant
under the action of the contraction $f$ can be written as
$$
x_i = \sum_{\cR'_i} a_{I,I'} (\vv,\overline{\vv})
$$
where each $a_{I,I'}:\CC^n \to \CC$ is a polynomial,
homogeneous of degree $I_k$ in $\vv_k$ and of degree
$I_k'$ in $\overline{\vv}_k$ for each $k$, such that
$a_{I,I'} + a_{I',I}$ is real-valued.
\end{thm}
\bigskip

\begin{rem}
(1) If $(J,J')$ belongs to $\cR'_i$ with real $\lambda_i$
above(!), then $(J', J)$ also belongs to $R_i'$. Consequently,
$a_{J,J'}+a_{J',J}$ is real-valued, as necessary, although
$a_{J,J'}$ and $a_{J',J}$ are in general complex-valued.
\medskip

(2) Most of the proof takes place at the formal (i.e., algebraic
and combinatorial) level; there is a small argument to take care
of a possible flat part, that we will present at the end of this
section.
\end{rem}
\bigskip

\subsection{Some preliminary observations}

We continue to use the special role of the index $i$ as above. Now
write
$$
\vv' = (\vv_1, \ldots, \vv_{i-1}; \vv'_i, y_i; \vv_{i+1}, \ldots,
\vv_m)
$$
so that $\vv = (x_i,\vv')$. The Taylor development of the defining
function of $M$ can be written as
$$
x_i = \sum a_{\alpha, I, I'} (\vv')
$$
where:
\begin{itemize} %
\item the sum is over multi-indices $(\alpha, I, I')$ such that
$\alpha + |I| + |I'| \ge 2$, (where $|I|$ denotes the sum of indices
constituting the multi-index $I$ for instance), and
\item $a_{\alpha,I,I'}$ is a real-valued polynomial homogeneous of
degree $\alpha$ in $y_i$, $I_k$ in $\vv_k$, $I'_k$ in
$\overline{\vv}_k$ $(k \neq i)$, $I_i$ in $\vv'_i$, and $I'_i$ in
$\overline {\vv'_i}$.
\end{itemize}
\bigskip

\begin{lem} \sl
If $a_{\alpha, I, I'} \neq 0$, then $\alpha = 0$, $I_k = I'_k = 0$
for every $k \ge i$.  In other words $x_i$, on the level of Taylor
series, depends only upon the variables $\vv_1, \ldots,
\vv_{i-1}$.
\end{lem}

{\it Proof}. Transforming $M$ by $f$, one obtains:
\begin{eqnarray*}
\widehat x_i & = & \lambda_i x_i + \hbox{Re } \sum_{I \in \cR_i}
c^0_I (\vv') \\
& = & \sum a_{\alpha,I,I'} (\widehat {\vv'}),
\end{eqnarray*}
together with
\begin{eqnarray*}
\widehat {\vv'} & = & (\widehat \vv_1, \ldots, \widehat \vv_{i-1},
\widehat {\vv'}_i, \widehat y_i, \widehat \vv_{i+1}, \ldots,
\widehat
\vv_m)  \\
\vspace{10pt} \\
\widehat \vv_j & = & \Lambda_j \vv_j + \sum_{J \in \cR_j} \cc_J (\vv) \\
\widehat y_i & = & \lambda_i y_i + \hbox{Im } \sum_{I \in \cR_i}
c_I^0 (\vv) \\
\widehat {\vv'}_j & = & \Lambda'_i \vv'_i + \cc' z_i + \sum_{I \in
\cR_i} \cc'_I (\vv) \qquad (j>i).
\end{eqnarray*}

Suppose, expecting a contradiction, that there exists $(\alpha, I,
I')$ with $a_{\alpha, I, I'} \neq 0$ and $\alpha + \sum_{k \ge i}
I_k + \sum_{k \ge i} I'_k > 0$. Amongst such terms
$a_{\alpha,I,I'}$, choose one which minimizes {\sl the total
degree of} $a_{\alpha,I,I'} = \alpha + \sum_\ell I_\ell +
\sum_\ell I'_\ell$. Then the invariance of $M$ under the action
by $f$ would imply that
$$
\lambda_i = \lambda_i^\alpha \lambda^I \overline{\lambda}^{I'}.
$$
But this last identity cannot hold, because we have
$$
|\lambda_j| \le |\lambda_i| \quad (\forall j \ge i),
$$
$$
\alpha + \sum_{k \ge i} I_k + \sum_{k \ge i} I'_k \ge 1,
$$
and
$$
\alpha + \sum_\ell I_\ell + \sum_\ell I'_\ell \ge 2.
$$
This contradiction yields the proof. \hfill $\Box$
\bigskip \\

We may now express $M$ (on the Taylor series level) by the
expression
$$
x_i = \sum a_{I,I'} (\vv_1, \ldots, \vv_{i-1}, \overline{\vv}_1,
\ldots, \overline{\vv}_{i-1})
$$
with $I_\ell = I'_\ell = 0$ for $\ell \ge i$. Notice that the
invariance under the action of $f$ yields
\begin{align}
\sum & a_{I,I'}(\vv_1,\ldots, \vv_{i-1}, \overline{\vv}_1, \ldots,
\overline{\vv}_{i-1}) \nonumber \\
& = \frac1{\lambda_i} \left( \sum a_{I,I'}
(\widehat{\vv}_1,\ldots, \widehat \vv_{i-1}, \overline{\widehat
\vv}_1, \ldots, \overline{\widehat
\vv}_{i-1}) \right) \label{equation-one} \\
& \qquad- \sum_{R_i} \hbox{Re } c_I^0 \nonumber
\end{align}

The main observation is the following

\begin{lem} \sl
If $a_{I,I'} (\vv_1, \ldots, \vv_{i-1}, \overline{\vv}_1, \ldots,
\overline{\vv}_{i-1})$ is homogeneous of degree $(I,I') \in
\cR'_i$ and if we set
$$
\widehat a_{I,I'} = a_{I,I'} (\widehat \vv_1, \ldots, \widehat
\vv_{i-1}, \overline{\widehat \vv}_1, \ldots, \overline{\widehat
\vv}_{i-1})
$$
(with $\displaystyle{\widehat \vv_j = \Lambda_j \vv_j +
\sum_{\cR_j} \cc_J (\vv)}$ as above), then the degree of each
non-zero homogeneous component of $\widehat a_{I,I'}$ also belongs
to $\cR'_i$.
\end{lem}
\bigskip

\it Proof. \rm Start with the relation
$$
\lambda_i = \lambda^I \overline{\lambda^{I'}}.
$$
When one develops $\widehat a_{I,I'}$ into homogeneous components,
each $\lambda_j$ or $\overline\lambda_j$ in the right hand term
either remains the same or is replaced by $\lambda^J$
(respectively $\overline\lambda^J$) for some $J \in \cR_j$. As
$\lambda_j = \lambda^J$ for $J \in \cR_j$ we still have equality.
Thus the corresponding degree belongs to $\cR'_i$. \hfill $\Box$
\bigskip \\

\subsection{Proof of Theorem \ref{theorem-main}}

We first prove that the Taylor series only contains components
with degree in $\cR'_i$.  If not, take the non-zero component with
degree not in $\cR'_i$ and with the smallest total degree.
Considering the corresponding components on both sides of
(\ref{equation-one}) gives a contradiction.
\medskip

It is now established that $M$ can be expressed by
$$
x_i = \sum_{\cR'_i} a_{I,I'} (\vv, \overline\vv) + \varphi(\vv')
$$
where $\varphi$ is $C^\infty$ and vanishes to infinite order at
$0$ (which is often called \it flat \rm at $0$).
\medskip

Assume that there exists $\vv'_{(0)}$ with $\varphi(\vv'_{(0)})
\neq 0$. Observe that, since $M$ is $f$-invariant, the
hypersurface
$$
M_0 ~:~ x_i = \sum_{\cR'_i} a_{I,I'} (\vv, \overline\vv)
$$
is also invariant.  Let $(x^0_{i,(0)}, \vv'_{(0)})$ (resp.,
$(x_{i,(0)}, \vv'_{(0))}$) be the point of $M_0$ (resp., $M$)
above $\vv'_{(0)}$. Write, for $\nu \ge 0$,
$$
f^\nu (x^0_{i,(0)}, \vv'_{(0)}) = (x^0_{i,(\nu)}, \vv'_{(\nu)}).
$$
The coordinate $x_{i,(\nu)}$ of $M$ that is above $\vv'_{(\nu)}$ is
represented by
$$
x_{i,(\nu)} = x^0_{i,(\nu)} + \varphi (\vv'_{(\nu)}).
$$
Denote its image under $f$ by $(\widetilde x_{i,(\nu+1)},
\widetilde \vv'_{(\nu+1)})$; see the drawing below:

\vspace{1cm}
 \setlength{\unitlength}{1mm}
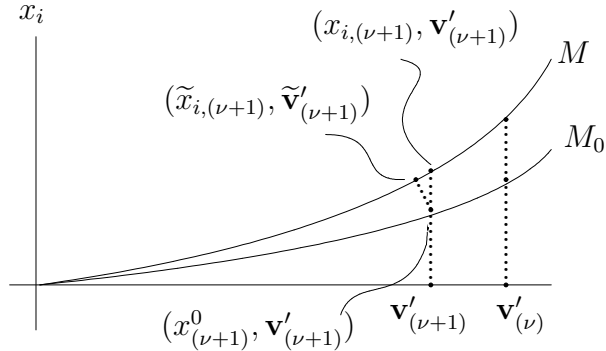
 \begin{center}
\begin{picture}(130,50)
\thinlines
\path(24.0,12.0)(24.0,12.0)(25.11,12.15)(26.22,12.31)(27.31,12.49)(28.4,12.65)(29.47,12.83)(30.56,13.0)(31.61,13.18)(32.66,13.36)
\path(32.66,13.36)(33.72,13.54)(34.75,13.72)(35.77,13.91)(36.79,14.11)(37.81,14.31)(38.81,14.5)(39.79,14.7)(40.79,14.9)(41.76,15.11)
\path(41.76,15.11)(42.73,15.33)(43.69,15.54)(44.62,15.75)(45.56,15.97)(46.51,16.18)(47.43,16.41)(48.33,16.63)(49.25,16.86)(50.13,17.09)
\path(50.13,17.09)(51.02,17.34)(51.91,17.56)(52.76,17.81)(53.62,18.06)(54.48,18.29)(55.33,18.54)(56.16,18.79)(56.98,19.04)(57.8,19.31)
\path(57.8,19.31)(58.61,19.56)(59.41,19.83)(60.19,20.09)(60.98,20.36)(61.75,20.62)(62.51,20.91)(63.26,21.18)(64.01,21.45)(64.76,21.75)
\path(64.76,21.75)(65.48,22.02)(66.19,22.31)(66.91,22.61)(67.62,22.9)(68.3,23.19)(68.98,23.48)(69.66,23.8)(70.33,24.09)(71.0,24.41)
\path(71.0,24.41)(71.63,24.72)(72.27,25.02)(72.91,25.34)(73.54,25.66)(74.15,25.98)(74.76,26.3)(75.34,26.62)(75.94,26.95)(76.51,27.3)
\path(76.51,27.3)(77.08,27.62)(77.65,27.97)(78.19,28.3)(78.75,28.65)(79.27,29.0)(79.8,29.34)(80.33,29.69)(80.83,30.05)(81.33,30.41)
\path(81.33,30.41)(81.83,30.76)(82.3,31.13)(82.77,31.5)(83.23,31.87)(83.69,32.23)(84.15,32.62)(84.58,32.98)(85.01,33.37)(85.43,33.75)
\path(85.43,33.75)(85.84,34.13)(86.25,34.52)(86.63,34.91)(87.02,35.3)(87.4,35.7)(87.76,36.11)(88.12,36.51)(88.48,36.91)(88.81,37.31)
\path(88.81,37.31)(89.15,37.73)(89.48,38.15)(89.79,38.55)(90.09,38.98)(90.4,39.41)(90.68,39.83)(90.95,40.26)(91.23,40.69)(91.5,41.12)
\path(91.5,41.12)(91.75,41.55)(91.98,41.98)(92.0,42.0)
\path(24.0,12.0)(24.0,12.0)(25.15,12.11)(26.29,12.24)(27.43,12.36)(28.56,12.47)(29.68,12.61)(30.77,12.74)(31.88,12.86)(32.97,12.99)
\path(32.97,12.99)(34.04,13.11)(35.11,13.25)(36.16,13.38)(37.22,13.52)(38.25,13.65)(39.29,13.79)(40.3,13.93)(41.33,14.06)(42.33,14.2)
\path(42.33,14.2)(43.31,14.34)(44.3,14.49)(45.27,14.63)(46.23,14.77)(47.19,14.93)(48.13,15.06)(49.06,15.22)(50.0,15.36)(50.91,15.52)
\path(50.91,15.52)(51.81,15.66)(52.7,15.83)(53.59,15.97)(54.48,16.13)(55.33,16.29)(56.19,16.45)(57.05,16.61)(57.88,16.77)(58.7,16.93)
\path(58.7,16.93)(59.52,17.09)(60.33,17.25)(61.13,17.41)(61.93,17.59)(62.7,17.75)(63.48,17.91)(64.25,18.09)(65.0,18.25)(65.73,18.43)
\path(65.73,18.43)(66.47,18.61)(67.19,18.77)(67.91,18.95)(68.62,19.13)(69.3,19.31)(69.98,19.5)(70.66,19.68)(71.33,19.86)(71.98,20.04)
\path(71.98,20.04)(72.63,20.22)(73.26,20.41)(73.9,20.59)(74.51,20.77)(75.12,20.97)(75.73,21.16)(76.3,21.34)(76.88,21.55)(77.45,21.73)
\path(77.45,21.73)(78.01,21.94)(78.56,22.12)(79.11,22.33)(79.65,22.52)(80.16,22.73)(80.68,22.93)(81.18,23.12)(81.66,23.33)(82.16,23.54)
\path(82.16,23.54)(82.62,23.75)(83.09,23.94)(83.55,24.16)(83.98,24.37)(84.43,24.58)(84.86,24.79)(85.26,25.01)(85.68,25.22)(86.06,25.43)
\path(86.06,25.43)(86.45,25.65)(86.83,25.87)(87.2,26.08)(87.56,26.3)(87.91,26.52)(88.25,26.75)(88.58,26.98)(88.9,27.19)(89.2,27.43)
\path(89.2,27.43)(89.51,27.65)(89.8,27.87)(90.08,28.11)(90.36,28.33)(90.62,28.58)(90.87,28.8)(91.12,29.04)(91.34,29.27)(91.58,29.51)
\path(91.58,29.51)(91.79,29.76)(91.98,29.98)(92.0,30.0)
\drawpath{20.0}{12.0}{92.0}{12.0} \drawpath{23.5}{45.0}{23.5}{6.0}
\drawcenteredtext{23.0}{48.0}{$x_i$}
\drawdotline{76.0}{12.0}{76.0}{27.2}
\drawdotline{86.0}{12.0}{86.0}{34.0}
\drawdotline{76.0}{22.0}{74.0}{26.0} \drawdot{74.0}{26.0}
\drawdot{76.0}{27.2} \drawdot{76.0}{22.0} \drawdot{76.0}{12.0}
\drawdot{86.0}{34.0} \drawdot{86.0}{26.0} \drawdot{86.0}{12.0}
\drawcenteredtext{95.0}{43.0}{$M$}
\drawcenteredtext{96.5}{31.5}{$M_0$}
\path(64.0,8.0)(64.0,8.0)(64.05,8.17)(64.11,8.34)(64.19,8.5)(64.25,8.65)(64.33,8.8)(64.41,8.94)(64.5,9.09)(64.58,9.21)
\path(64.58,9.21)(64.66,9.35)(64.75,9.46)(64.86,9.57)(64.94,9.69)(65.05,9.78)(65.15,9.88)(65.25,9.98)(65.36,10.07)(65.47,10.15)
\path(65.47,10.15)(65.58,10.23)(65.69,10.31)(65.81,10.39)(65.94,10.47)(66.05,10.52)(66.18,10.6)(66.3,10.65)(66.43,10.71)(66.55,10.77)
\path(66.55,10.77)(66.69,10.81)(66.81,10.86)(66.95,10.92)(67.08,10.96)(67.22,11.01)(67.36,11.05)(67.5,11.09)(67.63,11.13)(67.79,11.17)
\path(67.79,11.17)(67.93,11.21)(68.06,11.25)(68.22,11.27)(68.36,11.31)(68.5,11.35)(68.65,11.39)(68.8,11.43)(68.94,11.46)(69.09,11.5)
\path(69.09,11.5)(69.25,11.53)(69.4,11.57)(69.55,11.61)(69.69,11.65)(69.84,11.69)(70.0,11.73)(70.13,11.78)(70.29,11.84)(70.44,11.89)
\path(70.44,11.89)(70.58,11.93)(70.74,12.0)(70.88,12.05)(71.04,12.1)(71.19,12.18)(71.33,12.23)(71.47,12.31)(71.63,12.39)(71.77,12.47)
\path(71.77,12.47)(71.91,12.55)(72.05,12.63)(72.19,12.72)(72.34,12.81)(72.49,12.92)(72.63,13.02)(72.75,13.13)(72.9,13.23)(73.02,13.35)
\path(73.02,13.35)(73.16,13.48)(73.3,13.61)(73.43,13.76)(73.55,13.89)(73.68,14.05)(73.8,14.21)(73.93,14.36)(74.05,14.53)(74.16,14.71)
\path(74.16,14.71)(74.29,14.89)(74.4,15.09)(74.5,15.28)(74.61,15.48)(74.72,15.71)(74.83,15.93)(74.94,16.15)(75.04,16.39)(75.13,16.64)
\path(75.13,16.64)(75.22,16.89)(75.31,17.17)(75.41,17.35)(75.42,17.52)(75.43,17.82)(75.45,18.02)(75.47,18.33)(75.5,18.65)(75.52,18.99)
\path(75.54,19.29)(75.55,19.63)(75.60,19.38)(75.65,19.60)
\drawrighttext{66.0}{6.0}{$(x^0_{(\nu+1)}, {\vv}'_{(\nu+1)})$}
\drawrighttext{81.0}{8.0}{${\vv}'_{(\nu+1)}$}
\drawcenteredtext{88.0}{8.0}{${\vv}'_{(\nu)}$}
\path(62.0,32.0)(62.0,32.0)(62.22,31.98)(62.45,31.98)(62.68,31.97)(62.9,31.96)(63.11,31.95)(63.31,31.93)(63.5,31.9)(63.7,31.88)
\path(63.7,31.88)(63.88,31.86)(64.05,31.82)(64.22,31.79)(64.38,31.76)(64.55,31.71)(64.7,31.68)(64.86,31.62)(65.0,31.57)(65.13,31.53)
\path(65.13,31.53)(65.27,31.47)(65.41,31.43)(65.52,31.37)(65.65,31.3)(65.77,31.25)(65.88,31.19)(65.99,31.12)(66.08,31.05)(66.19,30.98)
\path(66.19,30.98)(66.27,30.92)(66.38,30.85)(66.45,30.77)(66.55,30.7)(66.63,30.62)(66.7,30.54)(66.79,30.46)(66.86,30.38)(66.93,30.29)
\path(66.93,30.29)(67.0,30.21)(67.05,30.13)(67.11,30.04)(67.18,29.96)(67.24,29.87)(67.3,29.79)(67.34,29.7)(67.41,29.62)(67.45,29.53)
\path(67.45,29.53)(67.5,29.44)(67.55,29.35)(67.61,29.26)(67.65,29.17)(67.69,29.07)(67.75,29.0)(67.79,28.9)(67.83,28.81)(67.88,28.72)
\path(67.88,28.72)(67.91,28.63)(67.97,28.54)(68.0,28.45)(68.05,28.37)(68.09,28.28)(68.13,28.19)(68.19,28.11)(68.22,28.02)(68.27,27.94)
\path(68.27,27.94)(68.33,27.85)(68.38,27.77)(68.43,27.69)(68.47,27.6)(68.52,27.52)(68.58,27.44)(68.63,27.37)(68.69,27.29)(68.75,27.21)
\path(68.75,27.21)(68.81,27.13)(68.88,27.06)(68.95,27.0)(69.02,26.93)(69.09,26.87)(69.16,26.79)(69.25,26.73)(69.33,26.68)(69.41,26.62)
\path(69.41,26.62)(69.5,26.55)(69.59,26.51)(69.69,26.45)(69.79,26.4)(69.9,26.36)(70.0,26.3)(70.11,26.27)(70.22,26.22)(70.34,26.2)
\path(70.34,26.2)(70.47,26.15)(70.59,26.12)(70.74,26.1)(70.86,26.07)(71.02,26.05)(71.16,26.04)(71.31,26.02)(71.47,26.01)(71.65,26.0)
\path(71.65,26.0)(71.81,26.0)(71.99,26.0)(72.5,26.0)
\drawcenteredtext{55.0}{36.0}{$(\widetilde x_{i,(\nu+1)},
\widetilde \vv'_{(\nu+1)})$ }
\path(72.0,42.0)(72.0,42.0)(72.22,41.93)(72.44,41.86)(72.66,41.81)(72.84,41.75)(73.04,41.68)(73.2,41.61)(73.36,41.54)(73.52,41.47)
\path(73.52,41.47)(73.65,41.4)(73.77,41.33)(73.9,41.25)(74.0,41.18)(74.09,41.11)(74.19,41.02)(74.27,40.95)(74.33,40.86)(74.38,40.79)
\path(74.38,40.79)(74.44,40.7)(74.49,40.61)(74.52,40.54)(74.55,40.45)(74.56,40.36)(74.58,40.27)(74.58,40.18)(74.58,40.09)(74.58,39.99)
\path(74.58,39.99)(74.56,39.9)(74.55,39.79)(74.52,39.7)(74.5,39.59)(74.45,39.5)(74.41,39.4)(74.38,39.29)(74.33,39.18)(74.27,39.06)
\path(74.27,39.06)(74.22,38.95)(74.16,38.84)(74.11,38.74)(74.05,38.61)(73.97,38.5)(73.91,38.38)(73.83,38.27)(73.77,38.15)(73.69,38.02)
\path(73.69,38.02)(73.61,37.9)(73.55,37.77)(73.47,37.63)(73.4,37.5)(73.31,37.38)(73.25,37.25)(73.16,37.11)(73.09,36.97)(73.02,36.83)
\path(73.02,36.83)(72.94,36.68)(72.88,36.54)(72.81,36.4)(72.75,36.25)(72.69,36.11)(72.63,35.95)(72.56,35.79)(72.52,35.65)(72.47,35.49)
\path(72.47,35.49)(72.41,35.33)(72.38,35.16)(72.33,35.0)(72.3,34.84)(72.27,34.68)(72.25,34.5)(72.24,34.34)(72.22,34.16)(72.22,33.99)
\path(72.22,33.99)(72.22,33.81)(72.24,33.63)(72.25,33.45)(72.27,33.27)(72.3,33.09)(72.33,32.9)(72.38,32.72)(72.44,32.52)(72.5,32.33)
\path(72.5,32.33)(72.58,32.13)(72.66,31.94)(72.75,31.73)(72.84,31.54)(72.95,31.32)(73.08,31.12)(73.2,30.92)(73.34,30.7)(73.5,30.48)
\path(73.5,30.48)(73.66,30.28)(73.83,30.05)(74.02,29.84)(74.22,29.62)(74.44,29.38)(74.66,29.17)(74.91,28.94)(75.16,28.7)(75.41,28.46)
\path(75.41,28.46)(75.69,28.22)(75.99,28.0)(76.0,28.0)
\drawcenteredtext{74.0}{46.0}{$(x_{i,(\nu+1)}, {\vv}'_{(\nu+1)})$}
\end{picture}
 \end{center}
 \vspace{1cm}

Here, the vector joining $(x^0_{i,(\nu+1)}, \vv'_{(\nu+1)})$ to
$(\widetilde x_{i,(\nu+1)}, \widetilde \vv'_{(\nu+1)})$ is (when
$\nu$ is sufficiently large) nearly vertical (i.e., nearly
parallel to the $x_i$-axis) while the vector joining $(\widetilde
x_{i,(\nu+1)}, \widetilde \vv'_{(\nu+1)})$ to $(x_{i,(\nu+1)},
\vv'_{(\nu+1)})$ is nearly horizontal. Therefore we have
$$
\lim_{\nu \to +\infty}
\frac{\varphi(\vv'_{(\nu+1)})}{\varphi(\vv'_{(\nu)})} = \lambda_i.
$$
On the other hand, the size of $\vv'_{(\nu)}$ approaches zero
exponentially fast; this is not compatible with the flatness of
$\varphi$ at $0$. Therefore, $\varphi$ has to be identically zero.
This completes the proof. \hfill $\Box$
\bigskip \\

\subsection{An Example} \label{eg}
Consider a holomorphic contraction $z \mapsto \widehat z:
(\CC^3,0) \to (\CC^3,0)$ with its linear part
$$
\left( \begin{array}{ccc}
  \lambda &  &  \\
    & \lambda^2 & \\
    &  &  \lambda^4
    \end{array} \right),
$$
for some real number $\lambda$ with $0 < \lambda < 1$. The
resonance sets are
$$
\cR_2 = \{(2,0,1)\} \hbox{ and } \cR_3 =\{ (4,0,1), (2,1,1), (0,2,1) \}.
$$
Hence in an appropriate coordinate system, $f$ is represented by
$$
(z_1, z_2, z_3) \mapsto (\widehat z_1, \widehat z_2, \widehat z_3)
$$
with
\begin{eqnarray*}
\widehat z_1 & = & \lambda z_1 \\
\widehat z_2 & = & \lambda^2 z_2 + Dz_1^2 \\
\widehat z_3 & = & \lambda^4 z_3 + A z_2^2 + B z_1^2 z_2 + C
z_1^4.
\end{eqnarray*}
The extended resonance set $\cR'_3$ is described as follows (Note here
that the extended resonances is written as $(a,a'; b,b'; \ldots)$
instead of the presentation $(a,b,\ldots; a',b',\ldots)$.):
\begin{multline*}
\cR'_3 = \big\{ (4,0; 0,0; 0,0), (0,4;0,0; 0,0), \\
(3,1;0,0; 0,0), (1,3;0,0; 0,0), \\
(2,2;0,0; 0,0), (2,0;1,0; 0,0), (0,2;0,1; 0,0), \\
(2,0;0,1; 0,0), (0,2;1,0; 0,0), \\ (1,1;1,0; 0,0), (1,1; 0,1; 0,0), \\
(0,0; 2,0; 0,0), (0,0; 0,2; 0,0), (0,0; 1,1; 0,0) \big\}
\end{multline*}
Now, we are looking only at the invariant hypersufaces tangent to the real
hyperplane represented by the equation $x_3=0$. The candidates for the
invariant hypersurface $M$ are consequently given by
\begin{eqnarray*}
z_3 + \overline{z}_3 & = & a z_1^4 + \bar a \bar z_1^4 + b z_1^3
\bar
z_1 + \bar b \bar z_1^3 z_1 + c z_1^2 \bar z_1^2 \\
& & + d z_1^2 z_2 + \bar d \bar z_1^2 \bar z_2 + e z_1^2 \bar z_2
+ \bar e \bar z_1^2 z_2 + f z_1 \bar z_1 z_2 + \bar f z_1 \bar z_1
\bar z_2 \\
& & + g z_2^2 + \bar g \bar z_2^2 + h z_2 \bar z_2
\end{eqnarray*}
where $c, h$ are real-valued.
\medskip

A direct calculation shows the following:
\medskip
\begin{itemize}
\item If such an example exists, then $A=0$.
\medskip
\item If $D=0$, then $A=B=C=0$. In particular, $f$ is linear, and
in this case any choices for $a, b, c, d, e, f, g, h$ will define
$M$ invariant under the action by the contraction.
\medskip
\item If $D \not= 0$, then the following relations should hold:
$$
f=h=0,~ \hbox{Re } e\bar D = 0,~ B = 2g\lambda^2 D,~ C=d\lambda^2 D +
D^2.
$$
In other words, $d, g, \lambda$ and $D$ determines $B$ and
$C$. $e$ just need to satisfy $\hbox{Re } e\bar D = 0$. $a, b \in
\CC$ and $c \in \RR$ can be arbitrarily chosen.
\medskip
\item If $D\not= 0$, then $M$ contains a non-trivial complex
analytic set. For instance with choices $\lambda=1/2, D=d=f=1,
e=\sqrt{-1}$, it follows that $B=1/2, C=5/4$. Then $M$ contains
the variety defined by $z_1=0, z_3 = z_2^2$. (Compare this with
the main theorem of \cite{Kim-Kim}.)
\end{itemize}

\section{CR structures admitting a CR contraction}

Recall that a CR structure (of hypersurface type) on a real $2n-1$
dimensional smooth manifold $M$ consists of a smooth hyperplane
distribution $\zeta = (\zeta_x)_{x \in M}$ (each $\zeta_x$ is of
dimension $2n-2$) in the tangent bundle $TM$ and an integrable
smooth almost complex structure $\widetilde J = (\widetilde
J_x)_{x \in M}$ on $\zeta$.
\medskip

Let $(\zeta,\widetilde J)$ be a CR structure in a neighborhood of
$0$ in $\RR^{2n-1}$.  Assume that there exists a smooth
contraction $\widetilde f$ at $0$ which preserves the CR
structure. We will see that after an appropriate smooth change of
coordinates, $\zeta, \widetilde J$ and $\widetilde f$ can be
written in a very simple form, so that in particular $M$ admits a
CR embedding into $\CC^n$.
\medskip

The first step is to use a theorem of Catlin (\cite{Catlin2}):
taking $\RR^{2n-1}$ as the hyperplane $\{\re z_n = 0\}$ in
$\CC^n$, there exists a smooth integrable almost complex structure
$J$ on $(\CC^n, 0)$ such that the CR structure on $(\RR^{2n-1})$
induced by $J$ coincides with $(\zeta, \widetilde J)$ up to a {\it
flat} error at $0$ (i.e., an error vanishing to infinite order at
$0$).
\smallskip

We actually need slightly more: we want that the contraction
$\widetilde f$ on $(\RR^{2n-1},0)$ can be extended to a
contraction $f$ of $(\CC^n, 0)$ which preserves $J$ up to an error
which is flat at $0$.
\medskip

The next step is to rectify $J$: there exists a smooth
diffeomorphism $\Phi_0$ of $(\CC^n, 0)$ such that ${\Phi_0}^* J$
is the {\sl standard} complex structure $J_{st}$ on $(\CC^n, 0)$.
Let $M_0 := \Phi_0^{-1} (\{\re z_n = 0\})$; it is equipped with a
CR structure $(\zeta_0, \widetilde J_0)$ that is the image of
$(\zeta, \widetilde J)$ under ${\Phi_0}^*$, which coincides,
up to a flat term, with the structure
$(\zeta_{M_0}, \widetilde J_{M_0})$ induced by $J_{st}$.

Let $f_0 := \Phi_0^{-1} \circ f \circ \Phi_0$. This is a local
smooth contraction of $(\CC^n, 0)$ that preserves $M_0$ and
$(\zeta_0, \widetilde J_0)$; it also preserves $J_{st}$, up to a
flat error at $0$.
\bigskip

We now use the normal form theorem (see \S
\ref{contraction-section}) for smooth contractions (on the finite
jet level only): given any $k \ge 1$, there exists a local
holomorphic diffeomorphism $\Phi_1$ of $(\CC^n, 0)$ such that the
local contraction $f_1 := \Phi_1^{-1} \circ f \circ \Phi_1$ has
Taylor expansion of order $k$ consisting only of resonant terms.
Take $k$ so large that all resonant terms have order strictly
smaller than $k-1$.

Let $\widehat f_1$ be the polynomial diffeomorphism of $\CC^n$ (of
degree $\ll k$) given by the Taylor expansion of $f_1$ of order
$k$. Let $M_1 = \Phi_1^{-1} (M_0)$ and let it be equipped with
$(\zeta_1, \widetilde J_1) := {\Phi_1}^* (\zeta_0, \widetilde
J_0)$. The hypersurface $M_1$ is invariant under $f_1$, and the
restriction of $f_1$ to $M_1$ preserves $(\zeta_1, \widetilde
J_1)$.
\medskip

Observe that the $k$-jet of $M_1$ at $0$ is therefore invariant
under $\widehat f_1$.  By the arguments on the hypersurfaces
invariant under holomorphic contractions in the proof of Theorem
\ref{theorem-main}, the $k$-jet of $M_1$ only contains the
resonant terms (of degree $\ll k$).  This $k$-jet defines a
polynomial hypersurface $\widehat M_1$ which is invariant under
$\widehat f_1$, and consequently also invariant under $f_1$ up to
a finite order $\ge k$.
\medskip

Let $(\widehat\zeta_1, \widehat J_1)$ be the canonical CR
structure on $\widehat M_1$; let $\Phi_2$ be a smooth
diffeomorphism of $(\CC^n, 0)$ which satisfies $\Phi_2 (M_1) =
\widehat M_1$ and coincides with the identity up to an order $\ge
k$ at $0$. Then ${\Phi_2}^* (\widehat\zeta_1, \widehat J_1) =
(\zeta_2, \widetilde J_2)$ is a CR structure on $M_1$ which
coincides with $(\zeta_1, \widetilde J_1)$ up to an order $\ge k$.
\medskip

It is invariant under $f_2 := \Phi_2^{-1} \circ \widehat f_1 \circ
\Phi_2$ (because $(\widehat \zeta_1, \widehat J_1)$ is invariant
under $\widehat f_1$).  Moreover, $f_1$ and $f_2$ coincide up to
an order $\ge k$.
\medskip

The following commutative diagram summarizes what we have done so
far:
\bigskip

$$
\begin{array}{ccccccccc}
\left({\RR^{2n-1} \atop \zeta, \widetilde J} \right) & {\Phi_0
\atop \longleftarrow} & \left({M_0 \atop \zeta_0, \widetilde J_0}
\right) & {\Phi_1 \atop \longleftarrow} & \left({M_1 \atop
\zeta_1, \widetilde J_1} \right) & {\Phi \atop \dashleftarrow} &
\left({M_1 \atop \zeta_2, \widetilde J_2} \right)  & {\Phi_2
\atop\longrightarrow} &
\left({\widehat M_1 \atop \widehat \zeta_1, \widehat J_1} \right)\\
&&&&&&&& \\
\downarrow \widetilde f &  &  \downarrow f_0 & &
\downarrow f_1 &
& \downarrow f_2 & & \downarrow \widehat f_1  \\
&&&&&&&& \\
\left({\RR^{2n-1} \atop \zeta, \widetilde J} \right) & {\Phi_0
\atop \longleftarrow} & \left({M_0 \atop \zeta_0, \widetilde J_0}
\right) & {\Phi_1 \atop \longleftarrow} & \left({M_1 \atop
\zeta_1, \widetilde J_1} \right) & {\Phi \atop \dashleftarrow} &
\left({M_1 \atop \zeta_2, \widetilde J_2} \right)  & {\Phi_2
\atop\longrightarrow} & \left({\widehat M_1 \atop \widehat
\zeta_1, \widehat J_1} \right)
\end{array}
$$
\bigskip

Here, $(\widehat M_1, \widehat \zeta_1, \widehat J_1 )$ is a
polynomial model. We now present:

\begin{lem}
There exists a local diffeomorphism $\Phi$ of $(M_1, 0)$ which
coincides with the identity up to terms of order $\ge k$ and
satisfies
$$
\Phi^* (\zeta_1, \widetilde J_1) = (\zeta_2, \widetilde J_2)
$$
and
$$
f_2 = \Phi^{-1} \circ f_1 \circ \Phi.
$$
\end{lem}

\it Proof. \rm The normal form theorem (cf.\ Proposition
\ref{contraction}, \S\ref{contraction-section}) for smooth
contractions guarantees that there exists a unique local
diffeomorphism $\Phi$ of $(M_1, 0)$ which coincides with the
identity up to terms of order $\ge k$, and satisfies $f_2 =
\Phi^{-1} \circ f_1 \circ \Phi$. Therefore in the lemma we can
assume that $f_2 = f_1$, and that $f_1$ (or, equivalently $f_2$)
preserves both $(\zeta_1, \widetilde J_1)$ and $(\zeta_2,
\widetilde J_2)$; these two CR structures coincide up to order
$\ge k$ and we have to show that they are in fact equal. \medskip

Let us first check that $\zeta_1 = \zeta_2$. Let $x$ be any point
close to $0$. We have
$$
\zeta_\ell (x) = [Df_1^N (x)]^{-1} (\zeta_\ell (f_1^N (x)))
$$
for $\ell = 1,2$. Let $\lambda < 1$ such that $\|f_1 (y)\| \le
\lambda \|y\|$ for $y$ close to $0$.  Then we have
$$
\| \zeta_1 (f_1^N(x)) - \zeta_2 (f_1^N (x)) \| \le C\lambda^{Nk}.
$$
As $k$ can be arbitrarily large, one can immediately conclude that
$\zeta_1 (x) = \zeta_2 (x)$.  The proof that $\widetilde J_1 =
\widetilde J_2$ is similar: if $v \in \zeta_1(x) = \zeta_2(x)$,
then one has
$$
\widetilde J_{\ell, x} (v) = [Df_1^N (x)]^{-1} (\widetilde
J_{1,f_1^N x} (Df_1^N (x) v))
$$
with
$$
\| \widetilde J_{1, f_1^N x} - \widetilde J_{2, f_1^N x} \| \le
C\lambda^{Nk}.
$$
Hence the assertion of lemma follows. \hfill $\Box\;$
\bigskip \\

As a consequence of the arguments by far, one obtains:

\begin{thm} \label{theorem-main-2}
Let $(M,0)$ be a germ of abstract smooth CR manifold of real
dimension $2n-1$ with the CR dimension $(2n-2)$. If there exists a
smooth CR contraction $f$ at $0$, then there exists a smooth CR
embedding $\psi:(M,0)\to (\CC^n,0)$ such that:
\begin{itemize}
\item[(1)] the image $\psi(M)$ is the hypersurface defined by a
real-valued weighted homogeneous polynomial;
\item[(2)] $\psi\circ f\circ \psi^{-1}$ coincides with the
polynomial normal form for the contraction $f$ at $0$ holomorphic
in a neighborhood of $0$ in $\CC^n$; and
\item[(3)] the degree of each monomial term in the expression of
the polynomial defining function for $\psi(M)$ is determined
completely by the extended resonance set of $f$.
\end{itemize}
\end{thm}

\section{Applications and Remarks}

\subsection{Strongly pseudoconvex case}
Assume that $(M,0)$ is a smooth, strongly pseudoconvex CR manifold
of hypersurface type, and that it admits a smooth CR contraction
$f$. Then Theorem \ref{theorem-main-2} implies that our $(M,0)$ is
CR equivalent to an embedded real hypersurface germ in $\CC^n$
that is invariant under a holomorphic contraction, and furthermore
that $(M,0)$ has to be CR equivalent to the strongly pseudoconvex
real hypersurface germ represented by
$$
\re z_n = Q(z_1, \ldots, z_{n-1})
$$
where $Q$ is a real-valued quadratic polynomial that is positive
definite (or, negative definite). Therefore, $(M,0)$ is CR
equivalent to the hypersurface represented by
$$
\re z_n = |z_1|^2 + \cdots + |z_n|^2.
$$
Notice that this yields an alternative proof of the
following theorem:

\begin{thm}[Schoen \cite{Schoen}]
Any smooth strongly pseudoconvex smooth CR manifold of
hypersurface type that admits a CR automorphism contracting at a
point is locally CR equivalent to the sphere.
\end{thm}

\subsection{Levi non-degenerate case}
In the case when the hypersurface germ $(M,0)$ is such that its
Levi form is not strongly pseudoconvex but only non-degenerate,
Schoen's theorem does not apply even when $M$ is assumed further
to be real analytic. For this particular case, Kim and Schmalz
have presented the following theorem:

\begin{thm}[\cite{Kim-Schmalz}]
If a real analytic, Levi non-degenerate hypersurface germ $(M,0)$
in $\CC^{n+1}$ ($n \ge 2$) is invariant under the action of a 1-1
holomorphic mapping that is repelling along the normal direction
to $M$ while contracting in some direction, then $(M,0)$ is
biholomorphic to the hyperquadric, defined by
$$
\re z_0 = |z_1|^2 + \ldots + |z_k|^2 - |z_{k+1}|^2 - \ldots -
|z_n|^2.
$$
\end{thm}

Then they give an example $\re z_3 = \re z_1 \overline{z}_2 +
|z_1|^4$ that defines a real analytic, Levi non-degenerate
hypersurface that is not biholomorphic to any hyperquadrics, while
it is invariant under the contraction $(z_1,z_2,z_3) \mapsto
(\lambda, \lambda^3, \lambda^4)$, in order to demonstrate that the
contraction is not relevant for their purposes.
\medskip

On the other hand, if one considers smooth ($\mathcal{C}^\infty$)
Levi non-degenerate CR manifold of hypersurface type with a CR
contraction, notice that this example of Kim-Schmalz {\sl does}
belong to the realm of Theorems \ref{theorem-main} and
\ref{theorem-main-2}.

Notice that, in the case that Levi form is indefinite (even when
it is non-degenerate) it is just that the Levi form alone cannot
determine the eigenvalues of the contracting CR automorphisms. But
if one starts from the given CR contraction, instead of Levi
geometric assumptions, Theorem \ref{theorem-main-2} reduces the
case to the embedded smooth CR hypersurface with a holomorphic
contraction, denoted by $f$ again, by an abuse of notation. Then
Theorem \ref{theorem-main-2} implies that the hypersurface germ
invariant under the action of $f$ is equivalent to the
hypersurface germ defined by an explicit polynomial defining
function with multi-degrees determined completely by the extended
resonance set for $f$. (See the example in \S\ref{eg} for
instance.)

\subsection{On the automorphism groups of bounded domains}
Let us now take a bounded domain $\Omega$ in $\CC^n$ with a smooth
($\mathcal{C}^\infty$) boundary $\partial\Omega$.  Assume further
the following:
\begin{itemize} \sl
\item[(1)]  $f$ is a holomorphic automorphism of $\Omega$.
\item[(2)] There exist a point $p \in \partial\Omega$ and an open
neighborhood $U$ of $p$ such that $f$ extends to a smooth
non-singular map of $U$ satisfying $f(p)=p$ and contracting at
$p$.
\end{itemize}
It follows that the iteration of $f$ in $U$ forms a sequence that
converges uniformly to the constant map with value at $p$. Then
Montel's theorem yields that the iteration of $f$ converges
uniformly on compact subsets of $\Omega$ to the same constant map.
\medskip

Choose a smooth extension of $f$ in a neighborhood $U$ of the
closure of $\Omega$. There is a smooth local diffeomorphism $\psi$
from a neighborhood $V$, say, of $p$ onto a neighborhood of $0$
such that $P := \psi \circ f \circ \psi^{-1}$ is a polynomial
normal form of $f$ at $0$.  Moreover, from Section 5.2.3, $\psi$
is holomorphic in $V \cap \Omega$. Let $U'$ be the basin of
attraction of $p$ for the extension of $f$; it is an open set
which contains $V$ and $\Omega$. One extends $\psi$ to a smooth
diffeomorphism (which we still denote by $\psi$) from $U'$ onto
its image by using the maps $P^{-k} \circ \psi \circ f^k$,
$k=1,2,\ldots$. Then $\psi$ still conjugates $f$ and $P$, and is
holomorphic in $\Omega$. The image under $\psi$ of the
intersection $U' \cap \partial\Omega$ is invariant under $P$, and
hence is a hypersurface of the form described in Theorem 2.3.
Consequently $\psi$ gives rise to a biholomorphism between
$\Omega$ and one side of this model hypersurface.

\medskip
Since $\Omega$ is bounded, one must have $i=n$. On the other hand,
this implies immediately that $\hbox{\rm Aut }(\widehat\Omega)$
admits free translation along $\hbox{\rm Im }z_n$-direction, as
well as the linear contractions. This yields the following:

\begin{cor}
Let $\Omega$ be a bounded domain with a smooth boundary. If there
exists $f \in \hbox{Aut }(\Omega) \cap \hbox{Diff
}(\overline\Omega)$ that is contracting at a boundary point, then
$\Omega$ is biholomorphic to the domain $\widehat\Omega$ defined
by a weighted homogeneous polynomial defining function. In
particular, the holomorphic automorphism group contains two
dimensional group of complex affine maps that are dilations
followed by translations.
\end{cor}

Notice that this is more general than the main results of
\cite{Kim-Kim}, in the sense that we do not put any restrictions
on the boundary except the smoothness. The first named author was
recently informed that a similar results as above (but on the
domains with D'Angelo finite type boundary) has been obtained by
S.-Y. Kim\cite{KIMSY} from a different argument, generalizing the
main result of \cite{Kim-Kim}.
\medskip

On the other hand our Theorem \ref{theorem-main-2} also yields the
following corollary:

\begin{cor}
Suppose that $\Omega$ is a bounded domain in $\CC^n$ with a smooth
boundary. If there exists $f \in \hbox{Aut }(\Omega) \cap
\hbox{Diff }(\overline\Omega)$ that is contracting at a boundary
point $p$, then $\partial\Omega$ at $p$ is of finite type in the
sense of D'Angelo \rm (\cite{DAN}).
\end{cor}

\it Proof. \rm Note that it suffices to show that the normal form
$\partial\widehat\Omega$ at the origin is of finite type in the
sense of D'Angelo. Since this boundary surface is defined by a
real-valued polynomial, notice further that it suffices to show
that there is no non-trivial complex analytic variety passing
through the origin. (See \cite{DAN2}.)
\medskip

Recall that $\partial\widehat\Omega$ is defined by the equation
$$
\re z_n = \rho(z_1, \ldots, z_{n-1})
$$
where $\rho$ is the real-valued polynomial consisting of its
monomial terms with degrees completely controlled by the set of
extended resonances of $f$. Consequently there exist positive
integers $1\le r_2 \le \cdots \le r_n$ such that
$\partial\widehat\Omega$ is invariant under the action of the
contraction
$$
P_\lambda :=
\begin{pmatrix}
 \lambda &               &         &       \\
         & \lambda^{r_2} &         &       \\
         &               & \ddots  &       \\
         &               &         & \lambda^{r_n}
\end{pmatrix}
$$
for every $\lambda$ with $0<\lambda<1$.
\medskip

Were there a non-trivial variety in $\partial\widehat\Omega$
passing through $0$, then one has a non-constant holomorphic map
$\varphi:D \to \partial\widehat\Omega$ of the open unit disc $D$
such that $\varphi(0)=0$.

Since $\partial\widehat\Omega \subset \CC^n$, we may write
$\varphi(\zeta) = (\varphi_1 (\zeta), \ldots, \varphi_n(\zeta))$.
For each $\mu \in \{1,\ldots,n\}$, consider the Taylor
developments of $\varphi_\mu (\zeta)$ at the origin. Let $d_\mu$
be the lowest degree term that does not vanish. (Set $d_\mu =
+\infty$ if $\varphi_\mu$ is identically zero.) Since $\varphi$ is
not identically zero, not every $d_\mu$ can be infinite.
\medskip

Now, consider the index $\widehat\ell \in \{1,\ldots,n\}$ such
that
$$
\frac{d_{\widehat\ell}}{r_{\widehat\ell}} = \min_{1 \le \ell \le
n} \frac{d_\ell}{r_\ell}.
$$

Let $R>1$ be arbitrarily given. Then let $\lambda =
R^{-{d_{\widehat\ell}}/{r_{\widehat\ell}}}$. Then we consider the
sequence of maps
$$
h_R (\zeta) := {P_\lambda}^{-1} \circ \varphi
\left(\frac{\zeta}{R}\right).
$$
It is a simple matter to observe that
$\displaystyle{\lim_{R\to\infty} h_R(\zeta)}$ defines a
non-constant entire curve with its image contained in
$\partial\widehat\Omega$.

But then, because the surface $\partial\widehat\Omega$ is defined
by the equation
$$
\re z_n = \rho(z_1, \ldots, z_{n-1}),
$$
one may add a real number to the last coordinate of the entire
curve generated above to obtain a non-constant entire curve
contained in $\widehat\Omega$. But this is again impossible
because $\widehat\Omega$ is biholomorphic to the bounded domain
$\Omega$. Thus the assertion follows. %
\hfill $\Box\;$
\bigskip

Notice that this gives an affirmative answer to a special case of
the Greene-Krantz conjecture in several complex variables, which
says that, for a bounded pseudoconvex domain with smooth boundary,
any boundary point at which an automorphism orbit accumulates is
of finite type in the sense of D'Angelo.
\medskip

As a passing remark, we would like to mention a naturally arising
question: \it when is the domain defined by such special
real-valued polynomial defining function biholomorphic to a
bounded domain? \rm Considering that it is relatively rare for a
basin of attraction to be bounded, this question may be of an
interest.

\subsection{A trivial remark}
It seems reasonable to add an obvious remark concerning the
point-wise attraction versus contraction. Even though the domain
$\Omega$ admits an automorphism $\varphi$ and a boundary point $p
\in \partial\Omega$ such that $\displaystyle{\lim_{\nu\to\infty}
\varphi^{\nu}(q) = p}$ for every $q \in \Omega$ (point-wise
attraction property) with the extra property that it has a smooth
extension to the boundary near $q$ fixing $q$, this extension may
not in general be a contraction at $p$. For each $t \in \RR$, the
holomorphic automorphism $z \mapsto t+z $ of the upper half plane
in $\CC$ with the unique fixed point at $\infty$. This gives rise
to an automorphism of the unit open ball in $\CC^2$ with
expression
$$
(z,w) \mapsto \left( c_t \cdot
\frac{z-\alpha_t}{1-\overline{\alpha_t} z},
\frac{\sqrt{1-|\alpha_t|^2}}{1-\overline{\alpha_t} z} w \right)
$$
where $c_t = (2i-t)/(2i+t)$ and $\alpha_t = t/(t-2i)$. Notice that
both eigenvalues of its derivative at the fixed point are $1$.
(This is a typical example of a parabolic map in complex dimension
2.)

\section{Appendix: Normal forms for smooth Contractions}

We recall some standard facts about smooth contractions, which
have been well-known for around one century.  The
exposition here follows that of \cite{PY}.
A simple fact about the holomorphicity of
the conjugacy to normal form is also proved at the end.

\subsection{Jets}

\subsubsection{}
Let $E$ be a real finite dimensional vector space. Denote by
$D(E,0)$ the group of germs of $\mathcal{C}^\infty$
diffeomorphisms of $E$ at $0$ which fix $0$. For $n\ge 0$, denote
by $D_n$ the normal subgroup of $D(E,0)$ whose elements have
contact of order $>n$ with the identity map $id_E$.  Namely we
have
$$
f(\xx) = \xx + o(\|\xx\|^n),
$$
for every $f \in D_n$.

Let $D_\infty := \bigcup_{n \ge 0} D_n$. For $0 \le n \le
+\infty$, the group of $n$-jets $J_n(E)$ is the quotient group
$D(E,0)/D_n$. It is obvious that $J_0 (E) = \{1\}$, $J_1 (E) =
GL(E)$. For $n \ge 1$, $J_{n+1}(E)$ is obtained from $J_n (E)$ by
an abelian extension
$$
1 \to D_n/D_{n+1} \to J_{n+1}(E)
\to J_n (E) \to 1,
$$
where $D_n/D_{n+1}$ is canonically identified with the polynomial
maps from $E$ to $E$ homogeneous of degree $n+1$. The groups
$J_n(E)$, $n<+\infty$, are Lie groups, and $J_\infty (E)$ is the
projective limit of the sequence $(J_n(E))_n$.

\subsubsection{} \label{1.2}
A jet $j$ in the Lie group $J_n (E)$ is semi-simple if and only if
it is conjugated to the $n$-jet of a semi-simple linear
automorphism of $E$; it is unipotent if and only if its image $J_1
(E)=GL(E)$ is a unipotent linear automorphism of $E$. Writing an
arbitrarily given jet as the product of its semi-simple and
unipotent components and going to the projective limit, we obtain
the following:
\medskip

\it
 For any $F \in D(E,0)$, one can find $H \in D(E,0)$, $F_1 \in
 D(E,0)$, $A \in GL(E)$ such that
 \begin{itemize}
 \item[(\romannumeral1)] $H \circ F \circ H^{-1} = F_1 \circ A$;
 \item[(\romannumeral2)] $A$ is semi-simple;
 \item[(\romannumeral3)] $DF_1(0)$ is unipotent;
 \item[(\romannumeral4)] $F_1 \circ A - A \circ F_1$ is flat at
 $0$.
 \end{itemize}
\rm
\bigskip

\subsubsection{} \label{1.3}
Let us describe the centralizer $Z(A)$ of a semi-simple linear
automorphism $A$ in the group $J_\infty(E)$, in the case where $A$
is a contraction.
\medskip

We first complexify the objects: let $E_\CC$, $A_\CC$, $J_\infty
(E_\CC)$ be these complexifications and let $Z_\CC (A)$ be the
centralizer of $A_\CC$ in $J_\infty (E_\CC)$.  Let $\lambda_1,
\ldots, \lambda_m$ be the {\sl distinct} eigenvalues of $A_\CC$,
and let $E_1, \ldots, E_m$ be the corresponding eigenspaces. As
$A_\CC$ is semi-simple, it follows that $E_\CC = \bigoplus E_\nu$.
For $\nu \in \{1,\ldots,m\}$ let
$$
R_\nu = \{I=(I_1, \ldots, I_m) \in \CC^m \mid
\sum I_\ell \ge 2, \lambda_\nu
= \lambda^I\}
$$
be the set of resonances. A jet $g \in J_\infty (E_\CC)$ with
$g(\vv_1,\ldots, \vv_m)= (g_1(\vv_1,\ldots, \vv_m),$ \newline $
\ldots, g_m(\vv_1,\ldots, \vv_m))$ belongs to $Z_\CC(A)$ if and
only if it can be written as
$$
g_\nu (\vv_1,\ldots, \vv_m) = B_\nu (\vv_\nu) + \sum_{I \in R_\nu}
g_{\nu, I} (\vv_1,\ldots, \vv_m)
$$
where:
\begin{itemize}
\item $B_\nu \in GL(E_\nu)$, and
\item $g_{\nu,I}$ is a polynomial map from $E$ to $E_\nu$
homogeneous of degree $I_\ell$ with respect to the vector variable
$\vv_\ell$ for $E_\ell$.
\end{itemize}
\medskip

Observe that, because $A$ is a contraction, the sets $R_\nu$ are
finite. It follows that $Z_\CC (A)$ is a finite dimensional Lie
group of polynomial diffeomorphisms of $E_\CC$.
\medskip

The real centralizer $Z(A)$ is then the real part of $Z_\CC (A)$,
i.e., formed by those jets which commute with complex conjugation.

\subsubsection{} \label{1.4}
There is a special case which we are interested in: assume that
$E$ is even-dimensional and equipped with a complex structure,
i.e., a linear automorphism $J$ such that $J^2 = -1_E$.  Assume
also that $A$ is $J$-linear, meaning that $AJ = JA$.  Then the
subgroup $Z_J (A)$ formed by those jets $g$ in $Z(A)$ which are
$J$-holomorphic has the form described above (where now
$\lambda_1, \ldots, \lambda_m$ are the eigenvalues of $A$
understood as a complex automorphism of $(E,J)$, $E_1, \ldots,
E_m$ being the $J$-invariant eigenspaces).

\subsection{The Normalization Theorem for Contractions}

\subsubsection{} \label{contraction-section}
Let $f \in D(E,0)$ be the germ of a smooth local contraction of
$E$ at $0$.  Write $A \in GL(E)$ for the semi-simple part of
$Df(0)$. According to \S\ref{1.2} above, there exists $h \in
D(E,0)$ such that the $\infty$-jet of $h^{-1} \circ f \circ h$
belongs to $Z(A)$.  On the other hand, by \S\ref{1.3}, $Z(A)$ is
canonically identified to a finite dimensional Lie group of
polynomial diffeomorphisms of $E$.

\begin{prop} \label{contraction}
The following statements hold:
\medskip

1. There exists $h \in D(E,0)$ such that $h^{-1} \circ f \circ h$
belongs to $Z(A)$ as an element of $D(E,0)$.
\medskip

2. Let $f_0, f_1 \in Z(A) \subset D(E,0)$.  Assume that both $f_0,
f_1$ have semi-simple part equal to $A$.  Then, any $h \in D(E,0)$
such that $h \circ f_0 \circ h^{-1} = f_1$ belongs to $Z(A)$.
\end{prop}

For the proof, see \cite{PY}.

\subsubsection{} \label{2.2}
We now make two simple observations in relation to Proposition
\ref{contraction}, one in this section and then the other in the
next.
\medskip

First, assume that $f_0, f_1$ are two local contractions in
$D(E,0)$ with the same $k$-th jet, where $k$ is large enough so that
$Z(A)$ injects into $J_k (E)$. (Here, $A$ is the semi-simple part
of $Df_0 (0) = Df_1 (0)$.) Then, there exists $h \in D(E,0)$ with
trivial $k$-th jet such that $h \circ f_0 \circ h^{-1} = f_1$.

\subsubsection{} The second observation, which we now give, is
more closely and explicitly related to the situation under
consideration of this paper.
\medskip

Assume as in \S\ref{1.4} that $E$ is even dimensional, and is
equipped with a complex structure $J$. Let $f \in D(E,0)$ be a
local contraction,and let $A$ the semi-simple part of $Df(0)$.
Assume also that the $\infty$-jet of $f$ is $J$-holomorphic (but
it is not necessarily a convergent series!).  Then there exists $h
\in D(E,0)$ with a $J$-holomorphic $\infty$-jet, such that $h
\circ f \circ h^{-1} \in Z_J (A)$.
\medskip

Assume even further that there exists an open set $U$ (with $0 \in
\overline{U}$) such that $f(U) \subset U$ and $f$ is
$J$-holomorphic on $U$. Then the conjugating map $h$ above is also
$J$-holomorphic on $U$. This can be see as follows: first
conjugate $f$ by a truncation $h_1$ of $h$ at a very high order
which is locally a biholomorphism at $0$.  We get a new
diffeomorphism $f_1 := h_1 \circ f \circ h_1^{-1}$ of the form
$$
f_1 = f_0 + o(\|\xx\|^k), k \gg 1
$$
where $f_0 \in Z_J (A)$. By the remark made in \S\ref{2.2} one has
that
$$
f_0 = h_0 \circ f_1 \circ h_0^{-1},
$$
$$
h_0 (\xx) = \xx + o(\|\xx\|^k).
$$
Now we are to check that $h_0$ is $J$-holomorphic on the open set
$U_1 := h_1(U)$.  Notice that $U_1$ is $f_1$-invariant, and on
$U_1$ the map $f_1$ is $J$-holomorphic. But for any $\xx \in U_1$
and $n \ge 0$, it holds that
$$
Dh_0 (\xx) = Df_0^{-n} (h_0 f_1^n (\xx))~ Dh_0 (f_1^n \xx)~ Df_1^n
(\xx).
$$
In this formula, both $Df_1^n(\xx)$ and $Df_0^{-n} (h_0 f_1^n
\xx)$ are $J$-linear and exponential in $n$.  On the other hand,
one has
$$
Dh_0 (f_1^n \xx) = 1_E + O(\|f_1^n x \|^k).
$$
Taking $k$ sufficiently large and letting $n$ go to infinity, one
deduces that $Dh_0(\xx)$ is $J$-linear.


\end{document}